\documentclass{article}[12pt]
\usepackage{amssymb}
\usepackage{amsfonts}
\usepackage{amsmath}
\usepackage[dvips]{graphics}
\input xy
\xyoption{all}
\newcommand{\vsp}{\vskip 1em}

\def \qed {\hfill \vrule height6pt width 6pt depth 0pt}

\newcommand{\bt}{\begin{theorem}}
\newcommand{\et}{\end{theorem}}
\newcommand{\be}{\begin{equation}}
\newcommand{\ee}{\end{equation}}
\newtheorem{theorem}{Theorem}[section]

\newtheorem{conjecture}[theorem]{Conjecture}
\newtheorem{corollary}[theorem]{Corollary}

\newtheorem{definition}[theorem]{Definition}
\newtheorem{example}[theorem]{Example}

\newtheorem{lemma}[theorem]{Lemma}

\newtheorem{problem}[theorem]{Problem}
\newtheorem{proposition}[theorem]{Proposition}

\begin{document}

\title{Maximizing spectral radius and number of spanning trees in bipartite  graphs}

\author{Ravindra B. Bapat\\
Indian Statistical Institute\\
New Delhi 110 016, India.\\
email: rbb@isid.ac.in
}

\date{}
\maketitle

\begin{abstract}
The problems of maximizing the spectral radius and the number
of spanning trees in a class of bipartite graphs with certain degree
constraints are considered. In both the problems, the optimal graph
is conjectured to be a Ferrers graph. Known results towards the resolution
of the conjectures are described. We give yet another proof of
a formula due to Ehrenborg and van Willigenburg for the number 
of spanning trees in a Ferrers graph. The main tool is a result
which gives several necessary and sufficient conditions under which
the removal of an edge in a graph does not affect the resistance distance
between the end-vertices of another edge. 
\end{abstract}


\vskip 1em


\newpage

\section{Introduction}

We consider simple graphs which have no loops or parallel edges.
Thus a { graph} $G = (V,E)$ consists of a finite set of
 vertices, $V(G),$ and a set of  edges, $E(G),$ each of whose
elements is a pair of distinct vertices.  
We will assume familiarity with basic graph-theoretic notions, see, for
example, Bondy and Murty [\ref{bonmurty}].

There are several matrices that one normally associates with a graph. We
introduce some  such matrices which are important.
Let $G$ be a graph with $V(G) = \{1,\ldots,n\}.$ 
The {\it adjacency matrix} $A$ 
of $G$ is an $n \times n$ matrix with  its rows and columns indexed by
$V(G)$ and with the $(i,j)$-entry equal to $1$ if vertices $i,j$ are
adjacent and $0$ otherwise. Thus $A$ is a
symmetric matrix with its $i$-th row (or column) sum equal to $d(i),$ which
by definition is the degree
of the vertex $i, i = 1,2,\ldots ,n.$ Let $D$ denote the $n \times n$
diagonal matrix, whose $i$-th diagonal entry is $d(i), i = 1,2,\ldots ,n.$
The {\it Laplacian matrix} of $G,$ denoted by $L,$ is  the matrix
$L = D - A.$

By the eigenvalues of a graph we mean the eigenvalues of its
adjacency matrix. Spectral graph theory is the study of the relationship
between the eigenvalues of a graph and its structural properties. The spectral
radius of a graph is the largest eigenvalue, in modulus, of the graph.
It is a topic of much investigation. It evolved 
during the study of molecular graphs by chemists. We refer to  [\ref{cvetkovic}] for
the subject of spectral graph theory.

A connected graph without a cycle is called a tree. Trees constitute
an important subclass of graphs both from  theoretical and  practical
considerations.
A spanning tree in a graph is a spanning subgraph which is a tree.
Spanning trees arise in several applications. If we are interested
in establishing a network of locations with minimal links, then 
it corresponds to a spanning tree. We may also be interested
in the spanning tree with the least weight, where each edge in
the graph is associated a weight and the weight of a spanning tree
is the sum of the weights of its edges.

If $G$ is connected, then $L$ is singular with rank $n-1.$ Furthermore, the well-known
Matrix-Tree Theorem asserts that any cofactor of $L$ equals the number of spanning
trees $\tau(G)$ in $G.$ For basic results concerning matrices associated with a graph
we refer to [\ref{bapat1}].

A graph $G$ is bipartite if its vertex set can be partitioned as $V(G) = X \cup Y$
such that no two vertices in $X,$ or in $Y,$ are adjacent. We often denote the bipartition
as $(X,Y).$ A graph is bipartite if and only
if it has no cycle of odd length.

The adjacency matrix of a bipartite graph $G$ has a particularly simple form viewed as a
partitioned matrix
$$A(G) = \left [ \begin{array}{cc}
0 & B \\ B' & 0 \end{array} \right ].$$

This form is especially useful in dealing with matrices associated
with a bipartite graph.

In this paper we consider two optimization problems over bipartite
graphs under certain constraints. One of the problems is to maximize
the spectral radius, while the other is to maximize the number
of spanning trees.

We now describe the contents of this paper. In Section 2 we introduce the
class of Ferrers graphs which are bipartite graphs such that the edges
of the graph are in direct correspondence
with the boxes in a Ferrers diagram. This class is of interest
in both the maximization problems that we consider.

The problem of maximizing the spectral radius of a bipartite graph
is considered in Section 3. We give a brief survey of the problem
and provide references to the literature containing results and open
problems.

In Section 4 we state an elegant formula for the number of spanning
trees in a Ferrers graph due to Ehrenborg and van Willigenburg [\ref{ehrenborg}].
We give references to the proofs of the formula available in the literature.
The formula leads to a conjectured upper bound for the number
of spanning trees in a bipartite graph and is considered in Section 5.
A reformulation of the conjecture in terms of majorization due to Slone
is described in Section 6.

Sections 7 and 8 contain new results. 
The concept of resistance distance [\ref{klein}] between two vertices in a graph
captures the notion of the degree of communication in a better way than
the classical distance. The resistance distance can be defined in several equivalent
ways, see, for example [\ref{bapat2}]. It is known, and intuitively
obvious, that the resistance distance between any two vertices 
does not decrease when an edge, which is not a cut-edge, 
is deleted from the graph.  In Section 7 we first give an introduction
to resistance distance. We then  examine the situation when
the removal of an edge in a graph does not affect the resistance distance
between the end-vertices of another edge. 
Several equivalent conditions are given for this to hold.
This result, which appears to be of interest by itself, is then
used in Section 8 to give another proof of the formula for the number
of spanning trees in a Ferrers graph.
Ehrenborg and van Willigenburg [\ref{ehrenborg}] also use electrical networks
and resistances in their proof of the formula but our approach is different.

\section{Ferrers graphs}

A Ferrers graph is defined as a bipartite graph on the bipartition $(U,V),$
where $U = \{u_1, \ldots, u_m\}, V = \{v_1, \ldots, v_n\}$ such that
\begin{itemize}
\item if $(u_i, v_j)$ is an edge, then so is $(u_p, v_q),$ where
$1 \le p \le i$ and $1 \le q \le j,$

\item $(u_1, v_n)$ and $(u_m, v_1)$ are  edges.

\end{itemize}

For a Ferrers graph $G$ we have the associated partition $\lambda = (\lambda_1, \ldots,
\lambda_m),$ where $\lambda_i$ is the degree of vertex $u_i, i = 1, \ldots, m,$
Similarly we have the dual partition $\lambda' = (\lambda_1', \ldots, \lambda_n')$
where $\lambda_j'$ is the degree of vertex $v_j, j = 1, \ldots, n.$ Note that
$\lambda_1 \ge \lambda_2 \ge \cdots \ge \lambda_m$ and $\lambda_1' \ge \lambda_2' \ge
\cdots \ge \lambda_n'.$ The associated
Ferrers diagram is the diagram of boxes where we have a box in position $(i,j)$
if and only if $(u_i,v_j)$ is an edge in the Ferrers graph.

\vskip 1em

\begin{example}  \label{example1}
The Ferrers graph with the degree sequences $(3,3,2,1)$ and $(4,3,2)$ is shown below.
\end{example}
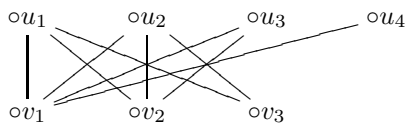
\[\xymatrix{\circ \mbox{$u_1$} \ar@{-}[d] \ar@{-}[rd] \ar@{-}[rrd] & \circ \mbox{$u_2$} 
\ar@{-}[ld] \ar@{-}[d] \ar@{-}[rd]& \circ \mbox{$u_3$} \ar@{-}[ld] \ar@{-}[lld] & \circ \mbox{$u_4$} \ar@{-}[llld] \\
 \circ \mbox{$v_1$} & \circ \mbox{$v_2$} & \circ \mbox{$v_3$}  
}\]
The associated Ferrers diagram is
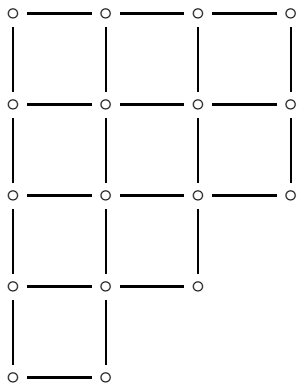
\[\xymatrix{\circ  \ar@{-}[d] \ar@{-}[r]  & \circ \ar@{-}[d] \ar@{-}[r] & \circ \ar@{-}[d] \ar@{-}[r] & \circ  \ar@{-}[d] \\
\circ  \ar@{-}[d] \ar@{-}[r]  & \circ \ar@{-}[d] \ar@{-}[r] & \circ \ar@{-}[d] \ar@{-}[r] & \circ  \ar@{-}[d] \\
\circ  \ar@{-}[d] \ar@{-}[r]  & \circ \ar@{-}[d] \ar@{-}[r] & \circ \ar@{-}[d] \ar@{-}[r] & \circ   \\
\circ   \ar@{-}[r] \ar@{-}[d] & \circ \ar@{-}[d] \ar@{-}[r] & \circ   \\
\circ   \ar@{-}[r] & \circ   
}\]

The definition of Ferrers graph is due to Ehrenborg and van Willigenburg [\ref{ehrenborg}].
Chestnut and Fishkind [\ref{chestnut}] defined  the class of bipartite graphs called
{\it difference graphs.}  A bipartite graph with parts X and Y is a difference graph if there
exists a function $\phi: X \cup Y \rightarrow R$ and a threshold $\alpha \in  R$
 such that for all $x \in X$ and $y \in Y,$ $x$ is adjacent to $y$  if and only if
$\phi (x) + \phi (y) \ge \alpha.$
It turns out that the class of Ferrers graphs coincides with the class of difference graphs,
as shown by Hammer et al. [\ref{hammer}].  A more direct proof of this 
equivalence is given by Cheng Wai Koo [\ref{koo}]. The same class is termed {\it chain 
graphs} in [\ref{amitava1}].

\section{Maximizing the spectral radius of a bipartite graph}

We introduce some notation.
Let $G = (V \cup W, E)$ be a bipartite graph,
 where $V = \{v_1, . . . , v_m\}, W = \{w_1, . . . ,w_n\}$ are the two partite sets. We
view the undirected edges $E$ of $G$ as a subset of $V \times W.$
Let $$D(G) = {d_1(G) \ge d_2(G) \ge \cdots \ge d_m(G)}$$ be the rearranged set of the degrees 
of $v_1, \ldots, v_m.$ Note that $e(G) = \displaystyle{\sum_{i=1}^m} d_i(G)$ 
is the number of edges in $G.$  Recall that the eigenvalues of $G$ are simply the eigenvalues
of the adjacency matrix of $G.$ Since the adjacency matrix is entrywise nonnegative, it follows from
the Perron-Frobenius Theorem that the spectral radius of the adjacency matrix is
an eigenvalue of the matrix. Denote by $\lambda_{max} (G)$ the maximum eigenvalue
of $G.$ 
It is known [\ref{amitava1}] that
\be \label{eq1.11}
\lambda_{max} (G) \le \sqrt{e(G)}
\ee
and equality occurs if and only if $G$ is a complete bipartite graph,
with possibly some isolated vertices.

We now consider refinements of (\ref{eq1.11}) for noncomplete bipartite graphs.
For positive integers $p,q,$ let $K_{p,q}$ be the complete bipartite
graph $G = (V \cup W, E)$ where $|V| = p, |W| = q.$
Let ${\cal K} (p, q, e)$ be  the family of subgraphs of $K_{p,q}$ with $e$ edges,  with no isolated
vertices, and which are not complete bipartite graphs.
The following problem was considered in [\ref{amitava1}]:

\begin{problem} \label{res31} Let $2 \le  p \le  q, 1 < e < pq$ be integers. Characterize the graphs which
solve the maximization problem
\be \label{eq11.2}
\max_{G \in  {\cal K}(p,q,e)} \lambda_{max}(G).
\ee
\end{problem}

Motivated by a conjecture of  Brualdi and Hoffman [\ref{brualdi}]  for nonbipartite
graphs, which was proved by Rowlinson [\ref{rowlinson}], the following conjecture
was proposed in [\ref{amitava1}]:

\begin{conjecture} \label{res32}
Under the assumptions of Problem \ref{res31},  an extremal graph that solves
the maximal problem (\ref{eq11.2})  is obtained from a complete bipartite graph by adding one vertex
and a corresponding number of edges.
\end{conjecture}

As an example, consider the class  ${\cal K}(3,4,10).$ 
There are two graphs in this class which satisfy the description
in Conjecture \ref{res32}. The graph $G_1$ obtained from the 
complete bipartite graph $K_{2,4}$ by adding an extra vertex of degree $2,$
and the graph $G_2,$ obtained from $K_{3,3}$ by adding an extra vertex of degree $1.$
The graph $G_1$ is  associated with the Ferrers diagram

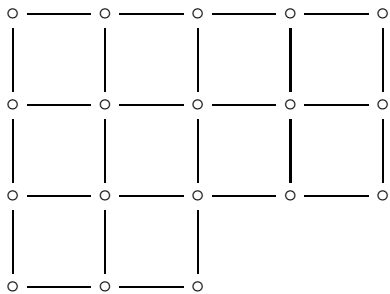
\[\xymatrix{\circ  \ar@{-}[d] \ar@{-}[r]  & \circ \ar@{-}[d] \ar@{-}[r] & \circ \ar@{-}[d] \ar@{-}[r] & \circ  \ar@{-}[d]
\ar@{-}[r] & \circ  \\
\circ  \ar@{-}[d] \ar@{-}[r]  & \circ \ar@{-}[d] \ar@{-}[r] & \circ \ar@{-}[d] \ar@{-}[r] & \circ  \ar@{-}[d] &
\circ \ar@{-}[l] \ar@{-}[u]\\
\circ  \ar@{-}[d] \ar@{-}[r]  & \circ \ar@{-}[d] \ar@{-}[r] & \circ \ar@{-}[d] \ar@{-}[r] & \circ   & \circ
\ar@{-}[u] \ar@{-}[l]\\
\circ   \ar@{-}[r]  & \circ  \ar@{-}[r] & \circ   
}\]
while $G_2$ is associated with the Ferrers diagram
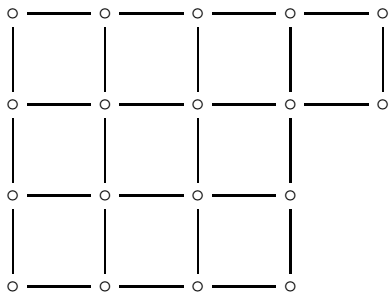
\[\xymatrix{\circ  \ar@{-}[d] \ar@{-}[r]  & \circ \ar@{-}[d] \ar@{-}[r] & \circ \ar@{-}[d] \ar@{-}[r] & \circ  \ar@{-}[d]
\ar@{-}[r] & \circ  \\
\circ  \ar@{-}[d] \ar@{-}[r]  & \circ \ar@{-}[d] \ar@{-}[r] & \circ \ar@{-}[d] \ar@{-}[r] & \circ  \ar@{-}[d] &
\circ \ar@{-}[l] \ar@{-}[u]\\
\circ  \ar@{-}[d] \ar@{-}[r]  & \circ \ar@{-}[d] \ar@{-}[r] & \circ \ar@{-}[d] \ar@{-}[r] & \circ   & 
\\
\circ   \ar@{-}[r]  & \circ  \ar@{-}[r] & \circ   & \circ \ar@{-}[u] \ar@{-}[l]
}\]

It can be checked that $\lambda_{max}(G_2) = 3.0592 > \lambda_{max}(G_1) = 3.0204.$
Thus according to Conjecture \ref{res32}, $G_2$ maximizes $\lambda_{max}(G)$ over
$G \in {\cal K}(3,4,10).$

Conjecture \ref{res32} is still open, although some special cases have been
settled, see [\ref{amitava1}, \ref{friedland}, \ref{petrovic}, \ref{stanley}]. We now mention a result from [\ref{amitava1}] 
toward the solution of Problem \ref{res31} which is of interest by itself, and is related to
Ferrers graphs.

Let $D = \{d_1, d_2, \ldots, d_m\}$  be a set of positive integers where
$d_1 \ge d_2 \ge \cdots \ge d_m$ and let ${\cal B}_D$ be the class of bipartite
graphs $G = (X \cup Y, E)$ with no isolated vertices, with $|X| = m,$ and with degrees of vertices
in $X$ being $d_1, \ldots, d_m.$ Then it is shown in [\ref{amitava1}] that 
$\max_{G \in {\cal B}_D} \lambda_{max} (G)$ is achieved, up to isomorphism, by the Ferrers graph,
with the Ferrers diagram having $d_1, d_2, \ldots, d_m$ boxes in rows $1,2, \ldots, m,$
respectively.

It follows that  an extremal graph
solving  Problem \ref{res31} is a Ferrers graph.

\section{The number of spanning trees in a Ferrers graph}

\begin{definition}  \label{res41} Let $G = (V, E)$ be a bipartite graph with bipartition $V = X \cup Y.$ The
Ferrers invariant of G is the quantity
$$F(G) = \frac{1}{|X||Y|} \prod_{v \in V} deg(v).$$
\end{definition}

Recall that we denote the number of spanning trees in a graph $G$ as $\tau(G).$
Ehrenborg and van Willigenburg  [\ref{ehrenborg}] proved the following interesting
formula.

\begin{theorem} \label{res41.1} If $G$ is a Ferrers graph, then $\tau(G) = F(G).$
\end{theorem}

Let $G$ be the Ferrers graph with  bipartition $(U,V),$ where $|U|=m, |V|=n.$ We assume
$U = \{u_1, \ldots, u_m\}, V = \{v_1, \ldots, v_n\}.$ 
Let $d_1 \ge \cdots \ge d_m$ and $d_1' \ge \cdots \ge d_n'$ be the degrees
of $u_1, \ldots, u_m$ and $v_1, \ldots, v_n$ respectively. We may assume $G$
to be connected, since otherwise, $\tau(G) = 0.$ If $G$ is connected, then
$d_1 = |V|$ and $d_1' = |U|.$ Thus according to Theorem \ref{res41.1},
$\tau(G) = d_2 \cdots d_m d_2' \cdots d_n'.$

As an example, the Ferrers graph in Example \ref{example1} has degree sequences
$(3,3,2,1)$ and $(4,3,2).$ Thus, according to Theorem \ref{res41.1}, it has
$3 \cdot 2 \cdot 1  \cdot 3 \cdot 2 = 36$ spanning trees.

 The complete graph $K_{m,n}$ has $m^{n-1}n^{m-1}$ spanning trees, and this can also
be seen as a consequence of Theorem \ref{res41.1}.

Theorem \ref{res41.1}  can be proved in many ways. The proof given by
Ehrenborg and van Willigenburg  [\ref{ehrenborg}] is based on electrical networks.
A purely bijective proof  is given by Burns [\ref{burns}]. We give yet another proof based on resistance
distance, which is different than the one in [\ref{ehrenborg}], see Section 8.

It is tempting to attempt a proof of Theorem \ref{res41.1} using the Matrix-Tree
Theorem. As an example, the Laplacian matrix of the Ferrers graph in Example
\ref{example1} is given by
$$L = \left [ \begin{array}{rrrrrrr}
3 & 0 & 0 & 0 & -1 & -1 & -1 \\
0 & 3 & 0 & 0 & -1 & -1 & -1\\
0 & 0 & 2 & 0 & -1 & -1 & 0\\
0 & 0 & 0 & 1 & -1 & 0 & 0\\
-1 & -1 & -1 & -1 & 4 & 0 & 0\\
-1 & -1 & -1 & 0 & 0 & 3 & 0\\
-1  & -1 & 0 & 0 & 0 & 0 & 2
\end{array} \right ].$$
Let $L(1|1)$ be the submatrix obtained from $L$ be deleting the first row and column. According to the 
Matrix-Tree Theorem, the number of spanning trees in the graph is equal to the determinant
of $L(1|1).$ Thus Theorem \ref{res41.1} will be proved if we can evaluate the determinant
of $L(1|1).$ But this does not seem easy in general.

A weighted analogue of Theorem \ref{res41.1} has also been given in [\ref{ehrenborg}] which
we describe now. 
Consider the  Ferrers graph $G$ on the vertex partition
$U = \{u_0, . . . , u_n\}$ and $V = \{v_0, . . . , v_m\}.$
For a spanning tree $T$ of $G,$ define the weight $\sigma(T)$ to be
$$\sigma (T) = \prod_{p=0}^n x_p^{deg_T(u_p)} \prod_{q=0}^m
y_q^{deg_T(v_q)},$$
where $x_0, \ldots, x_n; y_0, \ldots, y_m$ are indeterminates. 
 
For a Ferrers graph $G$ define $\Sigma(G)$ to be the sum $\Sigma (G) =
\sum_T \sigma (T),$  where $T$ ranges
over all spanning trees $T$ of $G.$ 

\begin{theorem} \label{res41.2} [\ref{ehrenborg}]
Let $G$ be the Ferrers graph corresponding to the partition $\lambda$ and the
dual partition $\lambda'.$
Then
$$\Sigma (G) =   x_0 \cdots  x_n \cdot y_0 \cdots  y_m
\prod_{p=1}^n (y_0 + \cdots ·+ y_{\lambda_p -1})
\prod_{q=1}^m (x_0 + \cdots +x_{\lambda'_q-1}).$$
\end{theorem}

Theorem \ref{res41.1} follows from Theorem \ref{res41.2} by setting
$x_0 = \cdots = x_n = y_0 = \cdots = y_m = 1.$

\section{Maximizing the number of spanning trees in a bipartite graph}

For general bipartite graphs the following conjecture was proposed 
by Ehrenborg [\ref{koo},\ref{slone}].

\begin{conjecture} \label{res42}(Ferrers bound conjecture). Let $G = (V, E)$ be a bipartite graph with
bipartition $V = X \cup Y.$ Then
$$\tau(G) \le \frac{1}{|X||Y|} \prod_{v \in V} deg(v),$$
that is, $\tau(G) \le F(G).$
\end{conjecture}

Conjecture \ref{res42} is open in general. In this section we describe
some partial results towards its solution, mainly from [\ref{garrett}]
and [\ref{koo}]. 
The following result has been proved in [\ref{garrett}].

\begin{theorem} \label{res50}
Let $G$ be a connected  bipartite graph for which Conjecture \ref{res42} holds. Let $u$ be a new
vertex not in $V(G),$ and let $v$ be a vertex in $V(G).$ Let $G'$ be the graph obtained by adding 
the edge $\{u,v\}$ to $G.$ Then Conjecture \ref{res42}  holds for $G'$ as well.
\end{theorem}

Note that Conjecture \ref{res42} clearly
holds for the graph consisting of a single edge. Any tree can be constructed
from such a graph by repeatedly adding a pendant vertex. Thus 
as an immediate consequence of Theorem \ref{res50} we get  the following.

\begin{corollary} \label{res51}
Conjecture \ref{res42}  holds when the graph  is a tree.
\end{corollary}

Using explicit calculations with homogeneous polynomials, the
following result is also established in [\ref{garrett}].

\begin{theorem} \label{res52} 
Let $G$ be a bipartite graph with bipartition $X \cup Y.$ Then
Conjecture \ref{res42} holds when $|X| \le 5.$
\end{theorem}

The following result is established in [\ref{koo}].

\begin{proposition} \label{res53}
Let $G$ and $G'$ be bipartite graphs for which Conjecture \ref{res42} holds.
Let $X$ and $Y$ be the parts of $G,$ and let $X'$ and $Y'$ be the parts of $G'.$
Choose vertices $x \in X$ and $x' \in X'.$ 
Define the graph $H$ with $V(H) = V(G) \cup V(G')$ and $E(H) = E(G) \cup E(G') \cup
\{xx'\}.$  Then the conjecture holds for $H$ also.
\end{proposition}

It may be remarked that Corollary \ref{res51} can be proved using Proposition
\ref{res53} and induction as well.
The following bound has been obtained in [\ref{bozkurt}].

\begin{theorem} \label{res54}
Let $G$ be a bipartite graph on $n \ge 2$ vertices. Then
\be \tau(G) \le \frac{\prod_v d_v}{|E(G)|}, \ee
with equality if and only if $G$ is complete bipartite.
\end{theorem}

Since there can be at most $|X||Y|$  edges in a bipartite graph with
parts $X$ and $Y,$  if Conjecture \ref{res42} were true, then Theorem 
\ref{res54} would follow.
Thus the assertion of Conjecture \ref{res42}  improves upon Theorem \ref{res54}
by a factor of $E(G)|/(|X||Y|).$ This motivates the following definition
introduced in [\ref{koo}].

\begin{definition}
 Let $G$ be a bipartite graph with parts $X$ and $Y.$ The bipartite density
of $G$, denoted $\rho (G),$ is the ratio $E(G)/(|X||Y|).$
Equivalently, $G$ contains $\rho (G)$ times as many edges as the complete
bipartite graph $K_{|X|,|Y|}.$  
\end{definition}

Let $G$ be a graph with $n$ vertices. Let $A$ be the adjacency matrix of $G$
and let $D$ be the diagonal matrix of vertex degrees of $G.$ Note that $L = D - A$
is the Laplacian of $G.$ The matrix $K = D^{-\frac{1}{2}}LD^{-\frac{1}{2}}$ is termed as
the normalized Laplacian of $G.$ If $G$ is connected, then $K$ is positive semidefinite
with rank $n-1.$ Let $\mu_1 \ge \mu_2  \cdots
\ge  \mu_{n-1} > \mu_n = 0$ denote the eigenvalues of $K.$ It is known, see [\ref{chung}],
that $\mu_{n-1} \le 2,$ with equality if and only if $G$ is bipartite. 
Conjecture \ref{res42} can be shown to be equivalent to the following, see [\ref{koo}].

\begin{conjecture} \label{res55}
 Let $G$ be a bipartite graph on $n \ge 3$ vertices with parts $X$ and
$Y.$ Then 
$$\prod_{i=1}^{n-2} \mu_i \le \rho (G).$$
\end{conjecture}

Yet another result from [\ref{koo}] is the following.

\begin{lemma}
Let $G$ be a bipartite graph on $n \ge 3$  vertices with parts $X$ and $Y.$
Suppose, for some $1 \le k \le \lfloor \frac{n-1}{2} \rfloor$ we have
$$\prod_{i=1}^k \mu_i(2 - \mu_i) \le \rho (G).$$
Then  Conjecture \ref{res42} holds for $G.$
\end{lemma}

We conclude this section by stating the following result [\ref{koo}].
It asserts that  Conjecture \ref{res42} holds for a sufficiently edge-dense graph
with a cut-vertex of degree $2.$

\begin{theorem}
 Let $G$ be a bipartite graph. Suppose that $\rho (G) \ge 0.544$ and that
$G$ contains a cut vertex $x$ of degree $2.$ Then Conjecture \ref{res42}  holds for $G.$
\end{theorem}

\section{A reformulation in terms of majorization} 

This section is based on [\ref{slone}].
Call a bipartite graph $G$ {\it Ferrers-good} if $\tau(G) \le F(G).$ Thus Conjecture \ref{res42}
may be expressed more briefly as the claim that all bipartite graphs are Ferrers-good.

In 2009, Jack Schmidt (as reported in [\ref{slone}]) computationally verified by an
exhaustive search that all bipartite graphs on at most $13$ vertices are Ferrers-good.
For a bipartite graph, we refer to the vertices in the two parts as red vertices
and blue vertices.
In 2013, Praveen Venkataramana proved an inequality weaker than Conjecture \ref{res42}
valid for all bipartite graphs:

\begin{proposition} \label{res43} (Venkataramana). Let G be a bipartite graph with red vertices having
degrees $d_1, \ldots, d_p$ and blue vertices having degrees $e_1, \ldots, e_q.$ Then
$$\tau(G) \le \prod_{i=1}^p (d_i + \frac{1}{2}) \prod_{j=1}^q (e_j + \frac{1}{2})
\sqrt{e_1}.$$
\end{proposition}

Conjecture \ref{res42} can be expressed in terms of majorization, for which the standard
 reference is [\ref{olkin}]. For a vector $a = (a_1, \ldots, a_n)$
the vector $(a[1], \ldots , a[n])$ denotes
the rearrangement of the entries of $a$ in nonincreasing order. Recall that a vector
$a = (a_1,\ldots, a_n)$ is majorized by another vector $b = (b_1, \ldots,b_n),$ written $a \prec b,$
provided that the inequality
$$\sum_{i=1}^k a_{[i]} \le \sum_{i=1}^k b_{[i]}$$
holds for $1 \le k \le n$ and holds with equality for $k = n.$

Given a finite sequence $a,$ let $\ell(a)$ denote its number of parts and $|a|$ denote its
sum. For example, if $a = (4, 3, 1),$ then $\ell(a) = 3$ and $|a| = 8.$

\begin{definition} \label{res45}  (Conjugate sequence). Let $a$ be a partition of an integer. The conjugate
partition of $a$ is the partition $a^*$ 
$$a_i^* = \# \{j : 1 \le j \le \ell (a) \mbox{ and } a_j \ge  i\}.
$$
\end{definition}

For example, $(5, 5, 4, 2, 2, 1)^* = (6, 5, 3, 3, 2).$

\begin{definition} \label{res46}  (Concatenation of sequences). Let $a = (a_1, \ldots, a_p)$ and 
$b = (b_1, \ldots,  b_q)$ be sequences. Then their conatentation is the sequence
$$a \oplus  b = (a_1, \ldots,  a_p, b_1, \ldots , b_q).$$
\end{definition}

With this notation, we can now state the following conjecture.

\begin{conjecture} \label{res47} Let $d$ be a partition with $\ell (d) = n,$ and let 
$\lambda$ be a non-increasing
sequence of positive real numbers with $\ell(\lambda) = n - 1.$ Suppose
$d = a \oplus b$ for some $a, b$ with $\ell (a) = p$ and $\ell (b) = q.$
If $a \prec b^*$ and $d \prec \lambda \prec d^*,$
then
$$\frac{1}{n} \prod_{i=1}^{n-1} \lambda_i \le \frac{1}{pq}
\prod_{i=1}^n d_i.$$
\end{conjecture}

Conjecture \ref{res47} implies  Conjecture \ref{res42} in view of the following two
theorems.

\begin{theorem} \label{res48} (Gale-Ryser). Let $a$ and $b$ be partitions of an integer. There is a
bipartite graph whose blue degree sequence is $a$ and whose red degree sequence is $b$ if
and only if $a \prec b^*.$
\end{theorem}

\begin{theorem} \label{res49}  (Grone-Merris conjecture, proved in [\ref{bai}])
The Laplacian spectrum of a
graph is majorized by the conjugate of its degree sequence.
\end{theorem}

Now let us show that  Conjecture \ref{res47} implies
Conjecture \ref{res42}. Assume Conjecture \ref{res47} is true.
 Let $G$ be a bipartite graph on $n$ vertices, with $p$ blue vertices and $q$ red
vertices. Let $d$ be its degree sequence, with blue degree sequence $a$ and red degree
sequence $b$, and let $\lambda$  be its Laplacian spectrum.
By Theorem \ref{res48}, $a \prec b^*.$  Since the Laplacian is a Hermitian matrix, 
$d \prec \lambda,$ and by Theorem \ref{res49}, $\lambda \prec  d^*.$
Hence the assumptions of Conjecture 2 apply.
We conclude that
\begin{equation} \label{eq41}
\frac{1}{n} \prod_{i=1}^{n-1} \lambda_i
\le \frac{1}{pq} \prod_{i=1}^n d_i.
\end{equation}
By the Matrix-Tree Theorem, the left-hand side of (\ref{eq41})  is $\tau(G).$ Hence
Conjecture \ref{res42}  holds as well.

\section{Resistance distance in $G$ and $G \setminus \{f\}$}

 We recall some definitions that will be useful.
Given a matrix $A$ of order $m \times n,$ a matrix $G$ of order $n \times m$
is called a generalized inverse (or a g-inverse) of $A$ if it satisfies $AGA = A.$
Furthermore $G$ is called Moore-Penrose inverse of $A$ if it satisfies
$AGA = A, GAG = G, (AG)' = AG$ and $(GA)' = GA.$ It is well-known that
the Moore-Penrose inverse exists and is unique. We denote the Moore-Penrose
inverse of $A$ by $A^+.$ We refer to [\ref{campbell}] for background material
on generalized inverses.

Let $G$ be a connected graph with vertex set $V = \{1, \ldots, n\}$ and let $i,j \in V.$
Let $H$ be a g-inverse of the Laplacian matrix $L$ of $G.$ The resistance distance
$r(i,j)$ between $i$ and $j$ is defined as
\be \label{eq1.1}
r_G(i,j) = h_{ii} + h_{jj} -h_{ij} - h_{ji}.
\ee
It can be shown that the resistance distance does not depend on the choice of the g-inverse.
In particular, choosing the Moore-Penrose inverse, we see that
$$r_G(i,j) = \ell^+_{ii} + \ell^+_{jj} - 2\ell^+_{ij}.$$

Let $G$ be a connected graph with $V(G) = \{1, \ldots, n\}.$
 We assume that each edge of $G$ is given an orientation.
If $e = \{i,j\}$ is an edge of $G$ oriented from $i$ to $j,$ then the incidence vector $x_e$
of $e$  is and $n \times 1$ vector with $1 (-1)$ at $i$-th ($j$-th) place and zeros elsewhere.
The Laplacian $L$ of $G$ has rank $n-1$ and any vector orthogonal to ${\bf 1}$
is in the column space of $L.$ In particular, $x_e$ is in the column space of $L.$

For a matrix $A,$ we denote by $A(i|j)$ the matrix obtained by deleting
row $i$ and column $j$ from $A.$ We denote $A(i|i)$ simply as $A(i).$
Similar notation applies to vectors.
Thus  for a vector $x,$ we denote by $x(i)$ the vector obtained by
deleting the $i$-th coordinate of $x.$
Let $L$ be the Laplacian matrix of a connected graph $G$ with vertex
set $\{1, \ldots, n\}.$  Fix
$i, j \in \{1, \ldots, n\}, i \not = j,$ and let $H$ be the matrix constructed as follows.
Set $H(i) = L(i)^{-1}$ and let the $i$-th row and column of $H$ be zero. Then
$H$ is a g-inverse of $L$ ([\ref{bapat1}], p.133). It follows from (\ref{eq1.1}) that
$r(i,j) = h_{jj}.$ For basic properties of resistance distance we refer to [\ref{bapat1},
\ref{bapat2}].

In the next result we give several equivalent conditions under which deletion of an edge 
does not affect the resistance distance between the end-vertices of another edge.
This result, which appears to be of interest by itself, will be used in Section 8
to give another proof of Theorem \ref{res41.1}. We denote an arbitrary g-inverse
of the matrix $L$ by $L^-.$

\begin{theorem} \label{res61} 
Let $G$ be a graph with $V(G) = \{1, \ldots, n\}, n \ge 4.$ Let
$e = \{i,j\}, f = \{k, \ell \}$ be  edges of $G$ with no common vertex such that $G \setminus \{e\}$ and
$G \setminus \{f\}$ are connected subgraphs. Let $L, L_e$ and $L_f$ be the  Laplacians of $G, 
G \setminus \{e\}$ and $G \setminus \{f\},$ respectively.
Let $x_e, x_f$ be the incidence vector of $e,f$ respectively. Then the 
following statements are equivalent:

\begin{description}

\item {(i)} $r_G(i,j) = r_{G \setminus \{f\}}(i,j)$

\item {(ii)} $r_G(k, \ell) = r_{G \setminus \{e\}} (k, \ell)$

\item {(iii)} $\tau (G \setminus \{e\}) \tau (G \setminus \{f\}) = \tau (G) \tau (G \setminus \{e,f\})$

\item {(iv)} The $i$-th and the $j$-th coordinates of $L^+x_f$ are equal

\item {(v)} The $i$-th and the $j$-th coordinates of $L^-x_f$ are equal for any $L^-$

\item {(vi)} The $i$-th and the $j$-th coordinates of $L_f^+x_f$ are equal

\item {(vii)} The $i$-th and the $j$-th coordinates of $L_f^-x_f$ are equal for any  $L_f^-$

\item {(viii)} The $k$-th and the $\ell$-th coordinates of $L^+x_e$ are equal

\item {(ix)} The $k$-th and the $\ell$-th coordinates of $L^-x_e$ are equal for any  $L^-$

\item {(x)} The $k$-th and the $\ell$-th coordinates of $L_e^+x_e$ are equal

\item {(xi)} The $k$-th and the $\ell$-th coordinates of $L_e^-x_e$ are equal for any $L_e^-.$
\end{description}
\end{theorem}

\noindent
{\bf Proof} Let $u = L_f^+x_f, w = L^+x_f.$  Since $x_f$ is in the column space of $L_f,$
we have $x_f = L_fz$ for some $z.$
It follows that $L_fu = L_fL_f^+x_f = L_fL_f^+L_fz
= L_fz = x_f.$
Similarly $Lw = x_f.$  Since $L = L_f + x_fx_f'$ then
$Lw = L_f w + x_fx_f'w$ and hence 
\be \label{eq101}
L_f(u-w) = x_fx_f'w.
\ee
Also, 
\be \label{eq102}
(x_f'w) L_fu = x_f(x_f'w).
\ee
Subtracting (\ref{eq102}) from (\ref{eq101}) gives
$L_f(u-w-(x_f'w))u = 0,$ which implies
$u-w - (x_f'w)u = \alpha {\bf 1}$ for some $\alpha.$ It follows that
$(1 - x_f'w)u = w + \alpha {\bf 1}.$ If $1 - x_f'w = 0,$ then all coordinates of $w$ are equal,
which would imply $Lw = 0,$ contradicting $x_f = Lw.$ Thus $1 - x_f'w \not = 0$ and hence
$u = \frac{w + \alpha {\bf 1}}{1 - x_f'w}.$ Thus any two coordinates of $u$ are equal if and only if the 
corresponding coordinates of $w$ are equal. This implies the equivalence of $(iv)$ and $(vi).$ 
A similar argument shows that $(iv)-(vii)$ are equivalent and that $(viii)-(xi)$ are equivalent.

Note that $r_G(i,j) = \frac {\det L(i,j)}{\det L(i)} = \frac{\tau (G \setminus \{e\})}{\tau (G)},
r_{G \setminus \{f\}} (i,j) = \frac {\det L_f(i,j)}{\det L_f(i)} = \frac{\tau (G \setminus \{e,f\})}{\tau (G\setminus \{f\})}.$
and
$r_{G \setminus \{e\}} (k, \ell) = \frac {\det L_e(k, \ell)}{\det L_e(k)} = \frac{\tau (G \setminus \{e,f\})}{\tau (G\setminus \{e\})}.$
Thus $(i),(ii)$ and $(iii)$ are equivalent.

We turn to the proof of $(iv) \Rightarrow (i).$ Let $w = L^+x_f$ and suppose 
$w_i = w_j.$ Since the vector ${\bf 1}$ is in the null space of $L^+,$ we may assume,
without loss of generality, that $w_i = w_j = 0.$
As seen before, $Lw = x_f.$

Since  $L(i) = L_f(i) + x_f(i)x_f(i)',$
by the Sherman-Morrison formula,
\begin{eqnarray}
L(i)^{-1} &=& (L_f(i) + x_f(i)x_f(i)')^{-1} \nonumber \\
&=& L_f(i)^{-1} - \frac{L_f(i)^{-1}x_f(i)x_f(i)'L_f(i)^{-1}} 
{1 - x_f(i)'L_f(i)^{-1}x_f(i)}. \label{eq1}
\end{eqnarray}
Since $x_f = Lw, w_i = 0$ and $(x_f(i))_j = 0,$  we have 
\begin{eqnarray*}
(x_f(i))_j &=& (L(i)w(i))_j \\
&=& ((L_f(i) + x_f(i)x_f(i)')w(i))_j\\
& =& (L_f(i)w(i)_j
+ x_f(i)'w(i)(L_f(i)x_f(i))_j.
\end{eqnarray*}
Hence $(L_f(i)^{-1}x_f(i))_j = 0.$ It follows from 
(\ref{eq1}) that the $(j,j)$-th element of $L(i)^{-1}$ and $L_f(i)^{-1}$
are identical.  In view of the observation preceding the Theorem,
the $(j,j)$-element of $L(i)^{-1}$ (respectively, $L_f(i)^{-1}$)
is the resistance distance between $i$ and $j$ in $G$ (respectively, $G \setminus \{f\}$).
Therefore the resistance distance between $i$ and $j$ is the same 
in $G$ and $G \setminus \{f\}$ if  the $i$-th and the $j$-th coordinates
of $L^+x$ are equal. 

Before proceeding we remark that if $(v)$ holds for a particular g-inverse, then it can be shown that
it holds for any g-inverse. Similar remark applies
to $(vii),(ix)$ and $(x).$

Now suppose $(i)$ holds. Then $(L(i))^{-1}_{jj} = (L_f(i))^{-1}_{jj},$
and using (\ref{eq1}) we conclude that
$(L_f(i)^{-1}x_f(i)x_f(i)'L_f(i))_{jj} = 0,$ which implies 
\be \label{eq2}
(L_f(i)^{-1}x_f(i))_{j} = 0.\ee

If we augment $L_f(i)^{-1}$ by introducing the $i$-th row and $i$-th column,
both equal to zero vectors, then we obtain a g-inverse $L_f^-$ of $L_f.$
Since the $i$-th coordinate of $x_f$ is zero, we conclude
from (\ref{eq2}) that $(L_f^-x_f)_{j} = 0.$ Since the $i$-th row of $L_f^-$
is zero, $(L_f^-x_f)_i = 0.$  It follows that the $i$-th and the $j$-th coordinates
of $L_f^-x_f = 0$ and thus $(vii)$ holds (for a particular g-inverse and hence for any g-inverse).
Similarly it can be shown that $(ii) \Rightarrow (xi).$ This completes the proof.
\qed

\section{The number of spanning trees in Ferrers graphs}

We now prove a preliminary result.

\begin{lemma} \label{lem1}
Consider the Ferrers graph $G$ with bipartition $(U,V),$ 
where $U = \{u_1, \ldots, u_m\}, V = \{v_1, \ldots, v_n\}.$
Let $\lambda_i$ be the degree of $u_i, i = 1, \ldots, m$ and let
$\lambda_j'$ be the degree of $v_j, j = 1, \ldots, n.$
Let $p \in \{1, \ldots, m-1\}$ be such that $\lambda_i = n, i = 1, \ldots, p$
and $\lambda_{p+1} = k < n.$ Let $f$ be the edge $\{u_p, v_n\}.$ Then
\be \label{eq10}
r_G(u_{p+1}, v_k) = r_{G \setminus \{f\}}(u_{p+1}, v_k).
\ee
\end{lemma}

\noindent
{\bf Proof} The  bipartite adjacency matrix of $G$ is given by
$$M = \bordermatrix{~ & 1 & 2 & \cdots & \cdots & n \cr
1 & 1 & 1 &  \cdots & \cdots & 1\cr
2 & 1 & 1 & \cdots & \cdots & 1 \cr
\vdots & 1 & 1 & \cdots & \cdots & 1  \cr
p & 1 & 1 & \cdots & \cdots & 1 & \cr
p+1 & 1 & 1 & \cdots & 0 & 0 \cr
\vdots & 1 & 1 & \cdots & 0 &  0 \cr
m & 1 & 1 & \cdots & \cdots & 0},$$
and the Laplacian matrix $L$ of $G$ is given by
$$L = diag(\lambda_1, \ldots, \lambda_m, \lambda_1', \ldots, \lambda_n') -
\left [ \begin{array}{cc}
0 & M \\ M' & 0 
\end{array} \right ].$$

Let
$$w = \frac{1}{p}[\underbrace{-\frac{1}{n}, \cdots, -\frac{1}{n}}_{p-1}, \frac{p-1}{n}, 0, \cdots, 0, -1]'.$$
It can be verified that  $Lw$ is the $(m+n) \times 1$ vector with $1$ at position $p,$ $-1$ at position $m+n$
and zeros elsewhere. Thus $Lw = x_f,$  the incidence vector of the edge $f = \{u_p,v_n\}.$ 

It follows from basic properties of the Moore-Penrose inverse [\ref{campbell}]
that 
$$ L^+L = \left (I - \frac{1}{m+n}{\bf 1}{\bf 1}' \right ).$$
Hence
\be \label{eq11}
L^+x_f = L^+Lw = \left (I - \frac{1}{m+n}{\bf 1}{\bf 1}' \right )w = w - \alpha {\bf 1}{\bf 1}',
\ee
where $\alpha = {\bf 1}'w/(m+n).$
Let $e$ be the edge $\{u_{p+1}, v_k\}.$ Since the coordinates $p+1$ and $m+k$ of $w$ are zero,
it follows from (\ref{eq11}) and the implication $(iv) \Rightarrow (i)$ of Theorem \ref{res61}
that (\ref{eq10}) holds. This completes the proof. \qed

\vskip 1em

Let $G$ be a connected graph with $V(G) = \{1, \ldots, n\},$ and let
$i,j \in V(G).$ Let $L$ be the Laplacian of $G.$ We denote by $L(i,j)$ the
submatrix of $L$ obtained by deleting rows $i,j$ and columns $i,j.$
Recall that  $\tau(G)$ denotes  the number of spanning trees of $G.$ It is well-known that
\be \label{eq2.1}
r_G(i,j) = \frac{\det L(i|j)}{\tau (G)}.
\ee
Furthermore, $\det L(i,j)$ is the number of spanning forests of $G$ with two
components, one containing $i$ and the other containing $j.$
Now suppose  that $i$ and $j$ are adjacent and let  $f = \{i,j\}$ be  the corresponding
edge. Let $\tau'(G)$ and $\tau''(G)$ denote the number of spanning trees of $G,$ containing
$f,$ and not containing $f,$ respectively. Then in view of the preceding remarks,
$\tau'(G) = \det L_1(i,j),$ where $L_1$ is the Laplacian of $G \setminus \{e\}.$

\bt  {[\ref{ehrenborg}]} Let $G$ be the  Ferrers graph with  the bipartition $(U,V),$
where $U = \{u_1, \ldots, u_m\}, V = \{v_1, \ldots, v_n\}$ and let 
$\lambda = (\lambda_1, \ldots,
\lambda_m), \lambda' = (\lambda_1', \ldots, \lambda_n')$
be the associated partitions. Then the number of spanning trees in $G$
is
$$\frac{1}{mn} \prod_{i=1}^m \lambda_i \prod_{i=1}^n \lambda_i'.$$
\et

\noindent
{\bf Proof} We assume $\lambda_m, \lambda_n'$ to be positive, for otherwise, the graph
is disconnected and the result is trivial. We prove the result by induction on the number of edges.
Let $e = \{p+1, m+k\}, f =\{p,m+n\}$ be edges of $G.$

By the induction assumption we have
\be \label{eq3}
\tau (G \setminus \{e\}) = \frac{1}{mn} \prod_{i=1}^m \lambda_i \prod_{i=1}^n \lambda_i' 
\frac{(\lambda_{p+1}-1)(\lambda_k' - 1)}{\lambda_{p+1} \lambda_k'},
\ee

\be \label{eq4}
\tau (G \setminus \{f\}) = \frac{1}{mn} \prod_{i=1}^m \lambda_i \prod_{i=1}^n \lambda_i' 
\frac{(\lambda_{p}-1)(\lambda_n' - 1)}{\lambda_{p} \lambda_n'},
\ee

and

\be \label{eq5}
\tau (G \setminus \{e,f\}) = \frac{1}{mn} \prod_{i=1}^m \lambda_i \prod_{i=1}^n \lambda_i' 
\frac{(\lambda_{p+1}-1)(\lambda_p)(\lambda_k'-1)(\lambda_n'-1)}{\lambda_{p+1} \lambda_p \lambda_k' \lambda_n'}.
\ee

It follows from (\ref{eq3}), (\ref{eq4}), (\ref{eq5}) and Theorem \ref{res61} that

$$
\tau (G) = \frac{\tau(G \setminus \{e\})(\tau (G \setminus \{f\})}{\tau (G \setminus \{e,f\})}
= \frac{1}{mn} \prod_{i=1}^m \lambda_i \prod_{i=1}^n \lambda_i',$$
and the proof is complete. \qed
 

\vsp

{\bf Acknowledgment} I sincerely thank Ranveer Singh for a careful reading of the manuscript.
Support from the JC Bose Fellowship, Department of Science and Technology, Government
of India, is gratefully acknowledged.

\vsp

\noindent
{\bf References}

\begin{enumerate}

\item \label{bai}
Hua Bai, The Grone-Merris conjecture, {\it  Transactions of the American 
Mathematical Society,} 363(8) (2011) , 4463–4474.

\item \label{bapat1}
R.B.  Bapat, {\it   Graphs and matrices},
Second edition, Hindustan Book Agency, New Delhi and  Springer,   2014. 

\item \label{bapat2}
R.B.  Bapat,   Resistance distance in graphs.
 Math. Student  68  (1999),  no. 1-4, 87--98.

\item \label{amitava1}
Amitava Bhattacharya, Shmuel friedland and Uri N. Peled,
On the first eigenvalue of bipartite graphs, {\it The Electronic Journal of Combinatorics,}
15 (2008) \# R144.

\item \label{bonmurty}
J.A.  Bondy and U.S.R.   Murty, U. S. R., {\it   Graph theory},
Graduate Texts in Mathematics, 244, Springer, New York,  2008.  

\item \label{bozkurt}
S.  Bozkurt and  Ş. Burcu,  Upper bounds for the number of spanning trees of graphs,
{\it J. Inequal. Appl.} 2012:269 (2012).

\item \label{brualdi}
R.A. Brualdi and A.J. Hoffman, On the spectral radius of (0, 1)-matrices, {\it Linear
Algebra Appl.}, 65 (1985), 133–146.

\item \label{burns} Jason Burns, Bijective proofs for "Enumerative Properties of Ferrers
Graphs", arXiv: math/0312282v1 [math CO] 15 Dec 2003.

\item \label{campbell} S.L. Campbell and C.D. Meyer, Jr., {\it
Generalized inverses of linear transformations}, Pitman, London, 1979.

\item \label{chestnut}
Stephen R. Chestnut and Donniell E. Fishkind,  Counting spanning
trees in threshold graphs, arXiv:1208.4125v2, (2013).

\item \label{chung}
R.K. Fan Chung, {\it  Spectral Graph Theory},  CBMS Regional Conference
Series in Mathematics, American Mathematical Society, 1997.

\item \label{cvetkovic}
D.M.  Cvetkovi\'{c}, Michael  Doob and Horst  Sachs, {\it Spectra of graphs.
Theory and applications}, Third edition, 
Johann Ambrosius Barth, Heidelberg,  1995.

\item \label{ehrenborg} Richard Ehrenborg and Stephanie  van Willigenburg, 
 Enumerative properties of Ferrers graphs, 
{\it  Discrete Comput. Geom.}  32  (2004),  no. 4, 481--492.

\item \label{friedland} S. Friedland, Bounds on the spectral radius of graphs with $e$ 
edges, {\it Linear Algebra Appl.} 101 (1988), 81–86.

\item \label{garrett} Fintan Garrett and Steven Klee,  Upper bounds for the
number of spanning trees in a bipartite graph. Preprint, \\
http://fac-staff.seattleu.edu/klees/web/bipartite.pdf, 2014.

\item \label{hammer}  Peter L. Hammer, Uri N. Peled, and Xiaorong Sun,
 Difference graphs, {\it  Discrete Applied Mathematics,} 28(1) (1990) 35 – 44.

\item \label{klein}
D.J. Klein and M.  Randi\'c,  Resistance distance, 
{\it  J. Math. Chem.},  12  (1993),  no. 1-4, 81--95.
	
\item \label{koo} Cheng Wai Koo,  A bound on the number of spanning trees in bipartite graphs,
2016. Senior thesis, https://www.math.hmc.edu/~ckoo/thesis/.

\item \label{olkin}
Albert W. Marshall, Ingram Olkin, and Barry C. Arnold. Inequalities: Theory
of Majorization and Its Applications. Springer, New York, 2011.

\item \label{rowlinson} 
P. Rowlinson, On the maximal index of graphs with a prescribed number of edges,
{\it Linear Algebra Appl.}, 110 (1988), 43–53.

\item \label{petrovic} 
Miroslav Petrovi\'c and  Slobodan K. Simi\'c,
A note on connected bipartite graphs of fixed order and size with maximal index,
{\it Linear Algebra Appl.},  483 (2015), 21-29.

\item \label{slone} Michael Slone, A conjectured bound on the spanning tree number of
bipartite graphs,  arXiv:1608.01929v2 [math.CO] 10 Aug 2016.
		
\item \label{stanley} R.P. Stanley, A bound on the spectral radius of graphs 
with $e$ edges, {\it Linear Algebra Appl.} 87 (1987), 267–269.

\end{enumerate}

\end{document}